\documentclass[12pt]{article}
\usepackage[landscape]{geometry}
\usepackage[francais]{babel}

\newtheorem{theoreme}{Th\'eor\`eme}
\newtheorem{proposition}{Proposition} 
\newtheorem{definition}{D\'efinition} 
\newtheorem{corollaire}{Corollaire}  
\newtheorem{notation}{Notation}
\newtheorem{lemme}{Lemme} 
 
\newtheorem{remarque}{Remarque}

 \title{L'id\'eal de Bernstein d'un arrangement libre d'hyperplans lin\'eaires}
 
 \author{ \vspace{1cm} \\ Ph. Maisonobe \\
   Universit\'e de Nice Sophia Antipolis \\
   Laboratoire Jean-Alexandre Dieudonn\'e \\
  Unit\'e Mixte de Recherche du  CNRS  7351 \\ 
  Parc Valrose, F-06108 Nice Cedex 2}

\begin{document}
\maketitle

 \newpage

\section*{Introduction}

 Soit $V$ un espace vectoriel de dimension $n$. Une famille $ \{ H_1, \ldots ,H_p \} $   d'hyperplans vectoriels    de $V$ deux \`a deux disctintcs   d\'efinit un arrangement
 ${\cal A} = {\cal A} ( H_1, \ldots ,H_p  ) $ de $V$.  Pour tout  $i \in \{1, \ldots ,p\}$,  soit $l_i$ une forme lin\'eaire sur $V$ de noyau $H_i$. Ces formes  sont d\'efinies
 \`a un coefficient de proportionalit\'e pr\`es.\\

 Suivons les notations du livre de P. Orlik et H. Terao \cite{O-T}. Notons $L( {\cal A} )$  l'ensemble des intersections non vide  d'\'el\'ements de  $ \{ H_1, \ldots ,H_p \} $.
  Si $X \in L( {\cal A} )$, nous notons   $ J(X) =\{ i \in \{1, \ldots ,p\} \; ; \; X \subset H_i\} $,   $r(x)$ le codimension de $X$ et
  ${\cal A}_X $ l'arrangement d'hyperplans d\'efini par $ \{ H_i \}_{ i \in J(X) } $. Un arrangement 
$ {\cal A}$ est dit irr\'eductible s'il n'est pas isomorphe \` a un produit d'arrangements   non vides.  
Soit $S = S(V^{\ast})$ l'alg\`ebre sym\'etrique de l'espace vectoriel dual de $V$.
Nous notons ${\rm Der}_{\bf C} (S)$ le $S$-module des d\'erivations de $S$ sur ${\bf C}$ et $ {\rm D} (  {\cal A} ) $ le  sous-module de ${\rm Der}_{\bf C} (S)$ form\'e des d\'erivations $\chi$ v\'erifiant
pour tout   $ i \in \{1, \ldots ,p\}$ :
$ \chi (f_i) \in l_i  \,  S  $.
 Nous disons   que l'arrangement  ${\cal A}$  est libre si $ {\rm D} (  {\cal A} ) $ est un $S$-module libre. \\
 
 En dimension $2$, tous les arrangements d'hyperplans sont libres. En dimension quelconque, les arrangements   
 d\'efinis par des sous-groupes de r\'eflexion du groupe lin\'eaire de $V$ d\'efinissent des arrangements libres dits  arrangements de reflexion. 
 C'est   le cas notamment des arrangements de tresse d\'efinis \` a partir des hyperplans invariants du groupe sym\'erique. \\
 
 Soit $A_V({\bf C})$, l'alg\`ebre de Weyl des op\'erateurs diff\'erentiels \`a coefficients dans $S$. Suivant la d\'emonstration de Bernstein \cite{Be}, l'id\'eal des polyn\^omes
 $ b\in {\bf C}[s_1, \ldots ,s_p]$ v\'erifiant :
 $$   \;\;   b (s_1, \ldots ,s_p) \,   l_1^{s_1} \ldots  l_p^{s_p} \in A_n({\bf C})[s_1, \ldots ,s_p] \,   l_1^{s_1+1} \ldots  l_p^{s_p+1} \; ,$$
 n'est pas  r\'eduit \`a z\'ero. Cet id\'eal  ne d\'epend pas du choix des formes lin\'eaires  $l_i$ qui d\'efinissent les hypersurfaces $H_i$. Nous notons cet id\'eal ${\cal B} ({\cal A} )$ et l'appelons l'id\'eal de Bernstein
 de ${\cal A} $. Le but de cet article est de d\'eterminer cet id\'eal  lorsque  ${\cal A}$  est un arrangement d'hyperplans lin\'eaires libre.\\
 
 Si $X  \in L(  {\cal A}  ) $, notons ${\cal A}_X $ l'arrangement d'hyperplans d\'efini par $ \{ H_i \}_{ i \in J(X) } $. 
 
 \newpage
 
 \begin{theoreme}
   Soit $V$ un espace vectoriel  et ${\cal A}$  un  arrangement d'hyperplans lin\'eaires   libre.   Alors  l'id\'eal  de Bernstein de    ${\cal A}$   est principal et engendr\'e par le polyn\^ome :
   $$ b_{\cal A} (s_1, \ldots ,s_p )= \prod_{X \in L'({\cal A })} \prod_{j=0}^{2( {\rm card } J(X) - r(X)) } ( \sum_{i \in J(X) } s_i  +r(X) +j   )   \quad ,   $$ 
   o\`u $L'({\cal A }) = \{ X  \in L( {\cal A }) \; ; \; {\cal A}_X \; {\rm est} \; {\rm irr\acute{e}ductible}  \}$. 
\end{theoreme}

\noindent Exemple 1 (en dimension 2) : Soit $l_1,l_2,\ldots ,l_p$ une famille de $p$   formes lin\'eaires homog\`enes sur ${\bf C}^2$ non deux \`a deux  colin\'eaires. L'id\'eal  de Bernstein de  $(  l_1, \ldots ,l_p)$   est principal et   engendr\'e par le polyn\^ome :
  $$ \prod_{i=1}^p   (s_i + 1)   \prod_{j=0}^{2(p-2)} ( s_1 + \cdots + s_p  + 2 + j ) \quad .$$
  
\noindent Exemple 2 (l'arrangement de Tresse) : Soit la famille $(x_i -x_j)_{1 \leq i <  j \leq n} $ de formes lin\'eaires sur ${\bf C}^n$.
L'id\'eal  de Bernstein de  $(x_i -x_j)_{1 \leq i <  j \leq n} $   est principal et   engendr\'e par le polyn\^ome :
  $$ \prod_{I \subset\{1,  \ldots ,n  \}}       \prod_{k=0}^{({\rm card}\, I -1)({\rm card} \, I -2)} ( \sum_{1 \leq i <  j \leq n \; ,\; i,j \in I} s_{i,j}  + {\rm card}\,  I -1 + k ) \quad .$$

 \section{Compl\'ements sur les champs de vecteurs logarithmiques}
 
 Soit $X$ une vari\'et\'e analytique complexe de dimension  $n$. D\'esignons par ${\cal O}_X$ son faisceau structural, 
  ${\cal D}_X$ le faisceau des op\'erateurs diff\'erentiels \` a coefficients dans  ${\cal O}_X$ et ${\rm Der}_X$ le sous-faiseau des champs de vecteurs holomorphes.\\

Soit $ (f_1 , \ldots , f_p)$  des fonctions analytiques sur $X$ et $F= f_1\ldots f_p$ leur produit. Notons $H_i$ l'hypersurface form\'ee par les z\'eros de la fonction $f_i$  et $H$ la r\'eunion de
ces hypersurfaces ou encore l'hypersurface form\'ee par les z\'eros de $F$. 
\\
 
 Soit  ${\cal D}_X [s_1,\ldots, s_p]$  le faisceau des op\'erateurs diff\'erentiels \` a coefficients dans  ${\cal O}_X [s_1,\ldots, s_p] $. D\'esignons par ${\cal D}_X [s_1,\ldots, s_p] f_1^{s_1}\ldots  f_p^{s_p}$, 
 le  ${\cal D}_X [s_1,\ldots, s_p]$-Module engendr\'e par $f_1^{s_1}\ldots  f_p^ {s_p}$
 dans ${\cal O}_X [s_1,\ldots, s_p, \frac{1}{F}] f_1^{s_1}\ldots  f_p^{s_p}$
 muni  de sa structure naturelle de 
 ${\cal D}_X [s_1,\ldots, s_p]$-Module.\\
 
 Notons par $T^{\ast}X \stackrel{\pi}{\rightarrow} X$  	le fibr\'e cotangent \` a $X$. Posons :
 $$\Omega_{f_1 , \ldots , f_p} =  \{ (x,  \sum_{i=1}^p s_i \frac{ df_i(x)}{f_i(x)}, s_1, \ldots ,s_p ) \; ; \; s_i \in {\bf C} \; ,  
\; {\rm et} \; F(x) \neq 0 \}\; .$$
 Cet ensemble $\Omega_{f_1 , \ldots , f_p} $  est sur   $\pi^{-1} (X-H ) \times {\bf C}^p$ une sous vari\'et\'e analytique lisse r\'eduite de dimension $n+p$ d\'efinie par les \'equations :
 $$ \xi_i - \sum_{i=1}^p s_i \frac{ df_i(x)}{f_i(x)} = 0  \; \; \; {\rm pour } \; i \in \{1, \ldots , p  \} \; .$$
 Nous
d\'esignons par $W^{\sharp}_{f_1 , \ldots f_p }$ l'adh\'erence dans $T^{\ast}X \times {\bf C}^p$ de  $\Omega_{f_1 , \ldots , f_p} $. Cette adh\'erence   est donc 
un sous-espace irr\'eductible de dimension  $n+p$ de  $T^{\ast}X \times {\bf C}^p$.\\

Suivant  K. Saito,  un champ de  vecteurs holomophe  $\chi $ est logarithmique pour $H$   si et seulement l'une des conditions \'equivalentes suivantes est v\'erifi\'ee :
\begin{itemize}
 \item  en tout point lisse de $H$, le  champ de  vecteurs  $\chi$ est tangent \` a $H$, 
\item pour tout $x_0$
dans  $H$, $\chi ( F)  \in F {\cal O}_{X,x_0}$.  
\end{itemize}
 Le sous-faisceau des d\'erivations logarithmiques pour $H$ est not\'e  ${\rm Der}_X (  H) $.
 De cette \'equivalence,  r\'esulte :
$$ {\rm Der}_X (  H)  = \cap_{i=}^p  {\rm Der}_X ( H_i )\; .     $$ 

\begin{notation} Soit $\chi \in {\rm Der}_X (  H)$ et    tout $j \in \{1, \ldots  , p  \}$, soit $b_j$ la fonction holomorphe
  telle que $ \chi (f_j) =b_j f_j$. Nous noterons $\tilde{\chi}$  l'op\'erateur diff\'erentiel  :
$$  \tilde{\chi} = \chi - \sum_{j=1}^p b_js_j  \in {\cal D}_X [s_1,\ldots, s_p]  $$
et par 
$$ \tilde{{\rm Der}}_X (f_1, \ldots,f_p)   = \{ \tilde{\chi} \; ; \; \chi \in {\rm Der}_X (   H)   \}    \; .$$
 
\end{notation}
Si $P$ est un op\'erateur de  ${\cal D}_X$, nous notons $\sigma (P)$ le symbole principal de $P$ pour la filtation naturelle de ${\cal D}_X$  par l'ordre des d\'erivations.
En donnant \`a chaque $s_i$ le poids un, cette filtration s'\'etend en une filtration de  ${\cal D}_X[s_1, \ldots ,s_p] $  que nous appelons  filtration di\`ese.
 Si $P \in {\cal D}_X[s_1, \ldots ,s_p] $,
   notons  $ \sigma ^{\sharp} (P)$   le symbole principal de $P$ pour cette filtration. C'est  une fonction analytique sur   $T^{\ast}X \times {\bf C}^p$.
   Si $M$ est un   ${\cal D}_X[s_1, \ldots ,s_p] $ -Module coh\'erent, nous noterons
  $ {\rm car}^{\sharp} (M)$ sa vari\'et\'e caract\'eristique relativement \`a la filration di\`ese.

    \begin{proposition} a) $ \tilde{{\rm Der}}_X (f_1, \ldots,f_p ) $ est l'ensemble des op\'erateurs de     ${\cal D}_X[s_1, \ldots ,s_p] $    d'ordre  un pour la filtration di\`ese annulant $ f_1^{s_1}\ldots  f_p^{s_p}$.\\
    \noindent b) L'application :
    $$  {\rm Der}_X (  H) \longrightarrow  \tilde{{\rm Der}}_X (f_1, \ldots,f_p) 
    \quad , \quad \chi \longmapsto \tilde{\chi}$$
    est un isomorphisme ${\cal O}_X$-lin\'eaire.
    \end{proposition}
    
    \begin{remarque} Si $U , V \in  {\rm Der}_X (   H)$ :
    $$ [U,V] = UV-VU \in {\rm Der}_X (H) \quad {\rm et} \quad   \tilde{ [U,V] } = \tilde{U} \tilde{V} - \tilde{V} \tilde{U} \; .$$
   \end{remarque}

    \begin{notation} \label{ngid} Nous notons $I^{\sharp}(f_1, \ldots ,f_p)$ l'id\'eal de
    ${\cal O}_{ T^{\ast}X }[s_1, \ldots ,s_p]$ engendr\'e par les symboles
    des op\'erateurs de $\tilde{{\rm Der}}_X (f_1, \ldots,f_p)$. Si $(\delta _1,
    \ldots ,\delta _l)$ forme un syst\`eme de g\'en\'erateurs du ${\cal O}_X$-Module
    ${\rm Der}_X (H) $ :
    $$I^{\sharp}(f_1, \ldots ,f_p)= ( \sigma^{\sharp} 
     (\tilde{\delta}_1),\ldots ,   \sigma^{\sharp}   (\tilde{\delta}_l) \quad .$$
        \end{notation}
        
     \begin{proposition} Si  $I^{\sharp}(f_1, \ldots ,f_p)$ est un id\'eal premier :
     \begin{enumerate} 
     \item la vari\'et\'e des z\'eros de $I^{\sharp}(f_1, \ldots ,f_p)$
     est $W^{\sharp}_{f_1 , \ldots f_p }\quad ,$  
      \item Nous avons (voir aussi \cite{B-M-M1} pour le r\'esutat g\'en\'eral) :
    $$  {\rm car}^{\sharp}{\cal D}_X [s_1,\ldots, s_p] f_1^{s_1}\ldots  f_p^{s_p}  = W^{\sharp}_{f_1 , \ldots f_p }   \quad ,  $$
     \item pour tout syst\'eme de g\'en\'erateurs $(\delta _1,
    \ldots ,\delta _l)$   du ${\cal O}_X$-Module
    ${\rm Der}_X (   H) $ :
    $${\rm Ann}_{ {\cal D}_X [s_1,\ldots, s_p]}  f_1^{s_1}\ldots  f_p^{s_p}
    = {\cal D}_X [s_1,\ldots, s_p] (   \tilde{\delta}_1 ,\ldots ,     \tilde{\delta}_l )   \quad .$$
    \end{enumerate}
   \end{proposition}
   
\noindent {\bf Preuve de 1:}    En dehors  de $H$, $ {\rm Der}_X (   H) $ est un ${\cal O}_{X-H}$-Module libre engendr\'e par les $n$  d\'erivations :
   $$ F \, ( \frac{\partial}{\partial x_i} - \sum_{j=1}^{p} \frac{1}{f_j}
   \frac{\partial f_j}{\partial  x_i} \, s_j ) \quad {\rm pour} \quad   1 \leq i \leq n   \quad ;$$
   Il en r\'esulte que la restriction de  $I^{\sharp}(f_1, \ldots ,f_p)$ 
   \`a $X-H$ est un id\'eal r\'eduit dont le lieu des z\'eros  est la vari\'et\'e lisse
    $\Omega_{f_1 , \ldots , f_p}$. Si   $I^{\sharp}(f_1, \ldots ,f_p)$ est un id\'eal premier, sa vari\'et\'e des z\'eros est irr\'eductible et est l'adh\'ernce de   $\Omega_{f_1 , \ldots , f_p}$ qui est par d\'efinition
    $W^{\sharp}_{f_1 , \ldots f_p }$.\\
 
   \noindent {\bf Preuve de 2:}   En fait ce r\'esultat est vrai sans hypoth\`ese sur $f_1,\ldots ,f_p$ voir \cite{B-M-M1}. Sous notre  hypoth\`ese, il se d\'eduit simplement du   premier r\'esultat.
   En effet un calcul direct montre que si un op\'erateur $P$ annule
   $f_1^{s_1}\ldots  f_p^{s_p}$, son symbole di\` ese est  nul sur 
   $\Omega_{f_1 , \ldots , f_p}$,  donc sur  $W^{\sharp}_{f_1 , \ldots f_p }$. Nous avons alors    :
   $$    W^{\sharp}_{f_1 , \ldots f_p } \subset   {\rm car}^{\sharp}{\cal D}_X [s_1,\ldots, s_p] f_1^{s_1}\ldots  f_p^{s_p} 
   \subset V(   I^{\sharp}(f_1, \ldots ,f_p) \subset  W^{\sharp}_{f_1 , \ldots f_p }  \quad .
   $$   
   
     \noindent {\bf Preuve de 3:}  Si $P \in {\cal D}_X [s_1,\ldots, s_p]$ annule  $  f_1^{s_1}\ldots  f_p^{s_p}$, $\sigma ^{\sharp}(P)$ appartient \'a l'id\'eal
     $( \sigma^{\sharp} 
     (\tilde{\delta}_1),\ldots ,   \sigma^{\sharp}   (\tilde{\delta}_l) )$. Ainsi, il existe $A_1, \ldots ,A_l \in {\cal D}_X [s_1,\ldots, s_p]$, tel que
     $$ \sigma^{\sharp}  (P) =  \sigma^{\sharp}  ( \sum_{j=1}^p  A_j \tilde{\delta}_j) \quad .
     $$ 
     L'op\'erateur$ P -    \sum_{j=1}^p  A_j \tilde{\delta}_j$ est un op\'erateur de degr\'e di\`ese strictement inf\'erieur \` a celui de $P$ annulant toujours
      $  f_1^{s_1}\ldots  f_p^{s_p}$. La d\'emonstration se fait alors par r\'ecurrence sur le degr\'e di\`ese de $P$.

\section {Le cas localement libre   quasi-homog\`ene}
Soit $X$ une vari\'et\'e analytique complexe de dimension  $n$,  
  $ (f_1 , \ldots , f_p)$  des fonctions analytiques sur $X$ et $F= f_1\ldots f_p$ leur produit.
 Nous conservons les notations du paragraphe pr\'ec\'edent. En particulier, $H_i$ d\'esigne l'hypersurface form\'ee par les z\'eros de la fonction $f_i$  et $H$ la r\'eunion de
ces hypersurfaces ou encore l'hypersurface form\'ee par les z\'eros de $F$.

 \begin{definition} $F$ est dit localement quasi-homog\`ene s'il existe au voisinage  de tout point un  syst\`eme de coordonn\'ees locales dans lequel $F= u \, G$
 o\`u $u$ est une unit\'e et $G$ est un polyn\^ome quasi-homog\`ene \`a poids stricitement positif. Le diviseur d\'efini  par $F$ est alors dit localement quasi-homog\`ene.
 \end{definition}

 \begin{proposition}  Si $F$ est localement quasi-homog\`ene  et que $ {\rm Der}_X (  H) $ est localement libre l'id\'eal
  $I^{\sharp}(f_1, \ldots ,f_p)$ est un id\'eal premier.
 \end{proposition}
 
      \noindent {\bf Preuve :} 
 Tout d'abors nous remarquons que si $u_1, \ldots, u_p$ sont des unit\'es de ${\cal O}_X$, $I^{\sharp}(f_1, \ldots ,f_p)$ est un id\'eal premier
 si et seulement si $I^{\sharp}(u_1f_1, \ldots ,u_pf_p)$ est un id\'eal premier.  Ainsi si  $x_0 \in H$, pour montrer notre  proposition au voisinage de $x_0$,
 nous pouvons supposer $X= {\bf C}^n$, $x_0=0$ et $F$ quasi-homog\`ene \`a poids strictement positif. Ainsi,    chaque $f_i$ sont quasi-homog\`ene \`a poids strictement positif.
 La preuve de la proposition se fait alors par r\'ecurrence.
 Le lemme clef, pour cette r\'ecurrence est la proposition de F. Castro, L. Narvaez et D. Mond :
 
 \begin{proposition}\cite{C-N-M} Soit $X$ une vari\'et\'e complexe de dimension $n$, soit $D$ un diviseur localement quasi-homog\`ene dans $X$ et $p \in D$.
 Alors, il existe un voisinage  ouvert $U$ de $p$ tel que pour tout $q \in U\cap D$, $ q \neq p$, le germe de paire $(X,D,q)$ est isomorphe au produit
 $({\bf C}^{n-1}\times {\bf C}, D' \times {\bf C}, (0,0))$ o\`u $D'$ est un diviseur localement quasi-homog\`ene. 
\end{proposition}

Supposons la proposition que nous voullons d\'emontrer vraie en dimmension $n-1$. 
Supposons $X = {\bf C}^{n}$ et $0 \in H$.  Il faut  montrer que le germe de l'id\'eal   $I^{\sharp}(f_1, \ldots ,f_p)$ est un id\'eal premier.
Soit $\delta _1, \ldots, \delta _n$ une base de $ {\rm Der}_X (   H) $, nous avons suivant la remarque \ref{ngid} :
$$  I^{\sharp}(f_1, \ldots ,f_p) = ( \sigma^{\sharp} 
     (\tilde{\delta}_1),\ldots ,   \sigma^{\sharp}   (\tilde{\delta}_n)) \quad .$$ Les composantes irr\'eductibles de la vari\'et\'e des z\'eros
     de $ I^{\sharp}(f_1, \ldots ,f_p)$ sont donc de dimension sup\'erieure \`a $n+p$. En dehors de $H$, $ I^{\sharp}(f_1, \ldots ,f_p)$
     cop\"{\i}ncide avce l'id\'eal premier :  
     $$  ( \xi _i - \sum_{j=1}^{p} \frac{1}{f_j}
   \frac{\partial f_j}{\partial  x_i} \, s_j ) \quad {\rm pour} \quad   1 \leq i \leq n   \quad ;$$
   Ainsi, $ W^{\sharp}_{f_1 , \ldots f_p }$ est une composante irr\'eductible de  la vari\'et\'e des z\'eros
     de $ I^{\sharp}(f_1, \ldots ,f_p)$. Si cette vari\'et\'e avait d'autres composantes irr\'eductibles, comme par r\'ecurrence
      $I^{\sharp}(f_1, \ldots ,f_p)$ est premier en dehors de $\pi^{-1}  (0)\times {\bf C}^{p}$, ses autres 
      composantes irr\'eductibles seraient contenues
     dans $\pi^{-1}(0)\times {\bf C}^{p}$. 
    Soit alors   $\chi$ le champ de quasi-homog\'enit\'e de $F$  et $d_i$ le degr\'e d'homog\'en\'eit\'e  de chaque $f_i$.
     Nous avons $\chi \in {\rm Der}_X (   H) $ et $ \tilde {\chi } = \chi  - \sum_{i=1}^p d_is_i$.
   Il en r\'esulte que  :
       $$( \pi^{-1}  (0)\times {\bf C}^{p} ) \cap  V (  I^{\sharp}(f_1, \ldots ,f_p) ) \;  \subset \;  T^{\ast}_{0}X \times (\sum_{i=1}^p d_i s_i =0) \quad . $$   
       Cette composante serait donc contenue dans un espace de dimension $n+p-1$. Nous en d\'eduisons :
       $$ V (  I^{\sharp}(f_1, \ldots ,f_p) )   =   W^{\sharp}_{f_1 , \ldots f_p } \quad .$$ 
       Il en r\'esulte que  $( \sigma^{\sharp} 
     (\tilde{\delta}_1),\ldots ,   \sigma^{\sharp}   (\tilde{\delta}_n))   $  est une suite r\'eguli\`ere. Comme $  I^{\sharp}(f_1, \ldots ,f_p) 
     $ est r\'eduit en dehors de $H$, il r\'esulte que  $ I^{\sharp}(f_1, \ldots ,f_p)  $ est r\'eduit et donc premier.

     \begin{corollaire} \label{crr} Si $F$ est localement quasi-homog\`ene  et que $ {\rm Der}_X (  \,  H) $ est localement libre, nous avons pour toute base
      $\delta _1, \ldots, \delta _n$   de $ {\rm Der}_X (   H) $  :
     \begin{enumerate}
     \item $I^{\sharp}(f_1, \ldots ,f_p)= ( \sigma^{\sharp} 
     (\tilde{\delta}_1),\ldots ,   \sigma^{\sharp}   (\tilde{\delta}_n)$ est un id\'esl premier et sa vari\'et\'e des z\'eros est $ W^{\sharp}_{f_1 , \ldots f_p }$,
     \item  $( \sigma^{\sharp} 
     (\tilde{\delta}_1),\ldots ,   \sigma^{\sharp}   (\tilde{\delta}_n))   $  est une suite r\'eguli\`ere de ${\cal O}_{T^{\ast}X}[s_1,\ldots, s_p]$,
          \item  ${ \sigma } 
     ( {\delta}_1),\ldots ,   \sigma   (  {\delta}_n))   $  est une suite r\'eguli\`ere de ${\cal O}_{T^{\ast}X} $,
    \item  $(F, \sigma^{\sharp} 
     (\tilde{\delta}_1),\ldots ,   \sigma^{\sharp}   (\tilde{\delta}_n))   $  est une suite r\'eguli\`ere de ${\cal O}_{T^{\ast}X}[s_1,\ldots, s_p]$, 
     \item $ {\rm Ann}_{ {\cal D}_X [s_1,\ldots, s_p]}  f_1^{s_1}\ldots  f_p^{s_p}
    = {\cal D}_X [s_1,\ldots, s_p] (   \tilde{\delta}_1 ,\ldots ,     \tilde{\delta}_l )  $,
    \item Soit $\dot{f_1}^{s_1}\ldots  \dot{f_p}^{s_p}$ la classe de $ f_1^{s_1}\ldots  f_p^{s_p}$ dans le quotient ${\cal D}_X [s_1,\ldots, s_p]   f_1^{s_1}\ldots  f_p^{s_p} /
    {\cal D}_X [s_1,\ldots, s_p]   f_1^{s_1+1}\ldots  f_p^{s_p+1}$, alors : 
  $ {\rm Ann}_{ {\cal D}_X [s_1,\ldots, s_p]}  \dot{f_1}^{s_1}\ldots  \dot{f_p}^{s_p}
    = {\cal D}_X [s_1,\ldots, s_p] ( F,  \tilde{\delta}_1 ,\ldots ,     \tilde{\delta}_n )  $.  
 \end{enumerate}
\end{corollaire}

   \noindent {\bf Preuve :} Les points $1$,$2$,$5$,$6$ r\'esulte directement ou facilement de ce qui pr\'ec\`ede. Pour le point $3$, rappelons suivant  
  par exemple \cite{B-M-M1} ou \cite{B-M-M2} que $ W^{\sharp}_{f_1 , \ldots f_p }\cap (s_1= \ldots =s_p=0)$ s'identifie \`a  sous-vari\'et\'e lagrangienen
  de $T^{\ast}X$. Ses composantes irr\'eductibles sont donc de dimension $n$. Il en r\'esulte que $(s_1, \ldots, s_p, \sigma^{\sharp} 
     (\tilde{\delta}_1),\ldots ,   \sigma^{\sharp}   (\tilde{\delta}_n)$ est un suite r\'eguli\`ere. Le r\'esultat  s'en d\'eduit. Pour le point 4,  cela r\'esulte du fait que par   d\'efinition
     de $ W^{\sharp}_{f_1 , \ldots f_p }$  :  si $u \in{\cal O}_X $ et $F\, u \in I^{\sharp}(f_1, \ldots ,f_p)$ alors  $u \in I^{\sharp}(f_1, \ldots ,f_p)$.\\

Dans la suite de ce paragraphe, nous supposerons $F$ est localement quasi-homog\`ene  et que $ {\rm Der}_X (   H) $ est localement libre. Notons 
${\rm Sp}^.(  \tilde{{\rm Der}}_X (f_1, \ldots,f_p)  )$  la suite de morphisme de   ${\cal D}_X [s_1,\ldots, s_p]$-Modules \`a gauche :
     $$ \longrightarrow {\cal D}_X [s_1,\ldots, s_p] \otimes_{{\cal O}_X[s_1,\ldots, s_p]} \Lambda  ^k ( \tilde{{\rm Der}}_X (f_1, \ldots,f_p)  ) \stackrel{ \epsilon _{-k}}{\longrightarrow}
{\cal D}_X [s_1,\ldots, s_p] \otimes_{{\cal O}_X[s_1,\ldots, s_p]} \Lambda  ^{k-1} ( \tilde{{\rm Der}}_X (f_1, \ldots,f_p)  ) \longrightarrow   $$
ou le morphisme $ \epsilon _{-k}$ est d\'efini par :
$$ \epsilon _{-k} ( P \otimes \eta _1 \wedge \ldots  \wedge \eta _k) = \sum_{k=1}^n    (-1)^i P \eta _i \otimes \eta _1 \wedge \ldots \wedge \hat{\eta _i} \wedge \ldots \wedge  \eta _k  
  + \sum_{1\leq i<j \leq k} (-1)^{i+j}  P   \otimes[ \eta _i , \eta _j ] \eta _1 \wedge \ldots \wedge \hat{\eta _i} \ldots \wedge \hat{\eta _j} \wedge \ldots \wedge  {\eta  _k}  \; .$$  
Suivant \cite{N}, ${\rm Sp}^.(  \tilde{{\rm Der}}_X (f_1, \ldots,f_p) )$ constitue 
  un complexe de  de   ${\cal D}_X [s_1,\ldots, s_p]$-Modules \`a gauche.

  \begin{proposition} \label{prldfs} Sous les hypoth\`eses $F$  localement quasi-homog\`ene  et   $ {\rm Der}_X (   H) $  localement libre, le complexe $ {\rm Sp}^.(   \tilde{{\rm Der}}_X (f_1, \ldots,f_p)  )$
  est une r\'esolution libre de ${\cal D}_X [s_1,\ldots, s_p]$-Module \`a gauche  du   Module  $   {\cal D}_X [s_1,\ldots, s_p] f_1^{s_1}\ldots  f_p^{s_p}$.
  \end{proposition}

     \noindent {\bf Preuve :} Ce complexe de ${\cal D}_X [s_1,\ldots, s_p]$-Module \`a gauche \` a gauche est naturelement filtr\'e . Son gradu\'e s'identifie au complexe de Koszul
      de la suite r\'eguli\`ere $( \sigma^{\sharp}  (\tilde{\delta}_1),\ldots ,   \sigma^{\sharp}   (\tilde{\delta}_n ))$. Il n'a donc de la cohomologie qu'en degr\'e z\'ero. Le morphisme d'augmention naturel
      (voir  le corollaire \ref{crr} ) :
     $$  {\rm Sp}^.(  \tilde{{\rm Der}}_X (f_1, \ldots,f_p)  )  \longrightarrow   {\cal D}_X [s_1,\ldots, s_p] f_1^{s_1}\ldots  f_p^{s_p}$$ 
     est donc un complexe  sans cohomologie.  
     
     \begin{notation} Soit $\chi \in {\rm Der}_X (   H)$ et pour    tout $i \in \{1, \ldots  , p  \}$, la fonction holomorphe
$b_i$ telle que $ \chi (f_i) =b_i f_i$. Nous noterons $\tilde{\tilde{\chi}}$  l'op\'erateur diff\'erentiel  :
$$ \tilde{\tilde{\chi}} = \chi + \sum_{j=1}^p b_j (s_j + 1)  \in {\cal D}_X [s_1,\ldots, s_p] \; ;$$
et par 
$$ \tilde{\tilde{{\rm Der}}}_X (f_1, \ldots,f_p)   = \{ \tilde{\tilde{\chi}} \; ; \; \chi \in {\rm Der}_X (   H)   \}    \; .$$
\end{notation}

Dans la suite,  nous supposserons que quitte \` a diminuer $X$,  $( \delta _1, \ldots , \delta _n )$ est une base de $ {\rm Der}_X (   H) $.
Il s'en suit que $( \tilde{\delta}_1 ,\ldots ,     \tilde{\delta}_n )  $
est une base du ${\cal O}]_X [s_1,\ldots, s_p] $-Module libre $ \tilde{{\rm Der}}_X (f_1, \ldots,f_p)$. De  m\^eme, $(  \tilde{\tilde{\delta}}_1 ,\ldots ,     \tilde{\tilde{\delta}}_n )  $
 est une base du ${\cal O}]_X [s_1,\ldots, s_p] $-Module libre $ \tilde{\tilde{{\rm Der}}}_X (f_1, \ldots,f_p)$.\\

   Soit $i$  et $j$ deux entiers distincts compris entre $1$ et $n$.
   Comme $ \delta _i  $ appartient \`a $ {\rm Der}_X (   H) $, il existe des $m_{i,j} \in {\cal O}_X$
   tels que :
   $$ \delta _ i (f_j) = m_{i,j} f_j\;.$$
   Nous ne d\'eduisons que 
  $\delta _i (F) = ( \sum_{j=1=}^n m_{i,j}) F$. De plus, 
   comme $ [\delta _i, \delta _j]$ appartient \`a $ {\rm Der}_X (   H) $, il existe des $\alpha _k^{i,j} \in {\cal O}_X$
  tels que :
  $$  [\delta _i, \delta _j] = \sum_{k=1}^n \alpha _k^{i,j} \delta _k \quad .$$ Nous en d\'eduisons :
  $$  [\tilde{\delta} _i, \tilde{\delta} _j] = \sum_{k=1}^n \alpha _k^{i,j} \tilde{\delta} _k 
  \quad {\rm et} \quad
   [\tilde{\tilde{\delta}} _i, \tilde{\tilde{\delta}} _j] = \sum_{k=1}^n \alpha _k^{i,j} \tilde{\tilde{\delta}} _k  \; .$$
   Notons que l'on : $\alpha _k^{i,j} = - \alpha _k^{j,i}$.
   
   \begin{lemme} Si  $^t \delta _i$ d\'esigne la transpos\'ee de l'op\'erateur $\delta _i$ :
   $$ ^t \tilde{\delta} _i  = - \tilde{\tilde{\delta}} _i + \sum_{k\neq i, k=1}^n \alpha _k^{i,k}\; .$$
   \end{lemme}

      \noindent {\bf Preuve :}  Cela repose sur un lemme \'etabli par F. Castro et J.-M. Ucha dans \cite{C-U} qui assure que pour tout $0\leq i\leq n$ :
   $$ ^t \delta _i =     - \delta _i - m_i + \sum_{k\neq i, k=1}^n \alpha _k^{i,k}\; ,$$
   Le lemem s'en d\'eduit puisque  $\tilde{\delta} _i  =   \delta _i - \sum_{j=1}^p m_{i,j} s_j$ et $\tilde{\tilde{\delta}} _i =  \delta _i + \sum_{j=1}^p m_{i,j} (s_j+1)$.\\

    Soit $M$ est un ${\cal D}_X [s_1,\ldots, s_p]$-Module coh\'erent. Notons $ \Omega ^n_X$ le faisceau des formes diff\'erentielles de degr\'e maximum sur $X$.
    Le dual de $M$ est le complexe $ RHom_{{\cal D}_X [s_1,\ldots, s_p]}(M , {\cal D}_X [s_1,\ldots, s_p])$  transform\'e en un complexe de
    ${\cal D}_X [s_1,\ldots, s_p]$-Module \` a gauche par application du foncteur $Hom_{{\cal D}_X [s_1,\ldots, s_p]}( \Omega ^n_X [s_1,\ldots, s_p] , - )$.
    Nous noterons $M^{\ast}$ ce complexe. Ses groupes de cohomologie sont les
    $$ Hom_{{\cal D}_X [s_1,\ldots, s_p]}( \Omega ^n_X [s_1,\ldots, s_p] , Ext^i_{{\cal D}_X [s_1,\ldots, s_p]}(M , {\cal D}_X [s_1,\ldots, s_p]))\; .$$

   \begin{proposition}  \label{pdual} Sous les hypoth\`eses $F$  localement quasi-homog\`ene  et   $ {\rm Der}_X (   H) $  localement libre :
   $$ ({\cal D}_X [s_1,\ldots, s_p] f_1^{s_1}\ldots  f_p^{s_p})^{\ast} \simeq {\cal D}_X [s_1,\ldots, s_p] f_1^{-s_1-1}\ldots  f_p^{-s_p-1} [-n] \; .$$
\end{proposition}

       \noindent {\bf Preuve :} Suivant \cite{M}, le nombre grade de ${\cal D}_X [s_1,\ldots, s_p] f_1^{s_1}\ldots  f_p^{s_p}$ est \'egal \`a $n$.
       Il r\'esulte alors de la propostion \ref{prldfs} que  le dual de $M$ est concentr\'e en degr\'e $n$. Calculons ce groupe de cohomologie.
       Explicitons pour cela :
  $$  {\cal D}_X [s_1,\ldots, s_p] \otimes_{{\cal O}_X[s_1,\ldots, s_p]} \Lambda  ^n ( \tilde{{\rm Der}}_X (f_1, \ldots,f_p)  ) \stackrel{ \epsilon _{-n}}{\longrightarrow}
{\cal D}_X [s_1,\ldots, s_p] \otimes_{{\cal O}_X[s_1,\ldots, s_p]} \Lambda  ^{n-1} ( \tilde{{\rm Der}}_X (f_1, \ldots,f_p)  )\; .$$
       Si $(\delta _1, \ldots, \delta _n)$ est une base de  $ {\rm Der}_X ({\rm Log} \,  H) $,  $\Lambda  ^n ( \tilde{{\rm Der}}_X (f_1, \ldots,f_p)  )$
       est un ${\cal O}_X [s_1,\ldots, s_p]$-Module libre  de base  $\tilde{\delta}_1 \wedge \ldots \wedge     \tilde{\delta}_n$ et
        $\Lambda  ^{n-1} ( \tilde{{\rm Der}}_X (f_1, \ldots,f_p)  )$ un ${\cal O}_X [s_1,\ldots, s_p]$-Module libre  de base
        $(\tilde{\delta}_1 \wedge \ldots \wedge \hat{\tilde{\delta}}_i  \wedge \ldots \wedge  \tilde{\delta}_n)_{1\leq i \leq n}$.
        Nous obtenons :
        $$ \epsilon _{-n} ( 1 \otimes \tilde{\delta}_1 \wedge \ldots \wedge     \tilde{\delta}_n ) 
        = \sum_{i=1}^n (-1)^{i-1} \tilde{\delta}_i \otimes \tilde{\delta}_1 \wedge \ldots \wedge \hat{\tilde{\delta}}_i  \wedge \ldots \wedge  \tilde{\delta}_n
        + \sum_{1\leq u <v \leq n } (-1)^{u+v} \otimes [\delta _u,  \delta _v] \wedge 
        \tilde{\delta}_1 \wedge \ldots \wedge \hat{\tilde{\delta}}_u   \wedge \ldots \wedge  \hat{\tilde{\delta}}_v  \wedge \ldots \wedge  \tilde{\delta}_n$$
          Et nous obtenons :
   $$ \epsilon _{-n} ( 1 \otimes \tilde{\delta}_1 \wedge \ldots \wedge     \tilde{\delta}_n ) 
        = (-1)^{i-1}  (  \tilde{\delta}_i - \sum_{k\neq i, k=1}^n \alpha _k^{i,k}) \otimes \tilde{\delta}_1 \wedge \ldots \wedge \hat{\tilde{\delta}}_i  \wedge \ldots \wedge  \tilde{\delta}_n\; .$$
       comme la transpos\'ee de l'op\'erateur $\tilde{\delta}_i - \sum_{k\neq i, k=1}^n \alpha _k^{i,k})$ est $ - \tilde{\tilde{\delta}} _i $. Nous obtenons que $M^{\ast}$
       est isomorphe au complexe concentr\'e en degr\'e $n$ :
       $$ \frac{ {\cal D}_X [s_1,\ldots, s_p]}{ {\cal D}_X [s_1,\ldots, s_p]  ( \tilde{\tilde{\delta}}_1 ,  \ldots ,   \tilde{\tilde{\delta}}_n )    } [-n] \; .$$
       Il reste \`a utiliser le point $4$ du corollaire \ref{crr}.
       
       \begin{remarque} Sous les hypoth\`eses de la proposition \ref{pdual}, nous pouvons en fait montrer l'existence d'un isomorphisme canonique entre 
       $$ Hom_{{\cal D}_X [s_1,\ldots, s_p]}( \Omega ^n_X [s_1,\ldots, s_p] ,  RHom_{{\cal D}_X [s_1,\ldots, s_p]}({\rm Sp}^.(  \tilde{{\rm Der}}_X (f_1, \ldots,f_p)  )  , {\cal D}_X [s_1,\ldots, s_p])  )\; .$$
       et
       $${\rm Sp}^.(  \tilde{\tilde{{\rm Der}}}_X (f_1, \ldots,f_p)  )[-n]  \; .$$

       \end{remarque}
       
          \begin{proposition} \label{pdq} Sous les hypoth\`eses $F$  localement quasi-homog\`ene  et   $ {\rm Der}_X (   H) $  localement libre :
   $$ ( \frac{ {\cal D}_X [s_1,\ldots, s_p] f_1^{s_1}\ldots  f_p^{s_p}}{ {\cal D}_X [s_1,\ldots, s_p] f_1^{s_1+1}\ldots  f_p^{s_p+1}})^{\ast}
   \simeq\frac{ {\cal D}_X [s_1,\ldots, s_p] f_1^{-s_1-2}\ldots  f_p^{-s_p-2}}{ {\cal D}_X [s_1,\ldots, s_p] f_1^{-s_1-1}\ldots  f_p^{-s_p-1}} [-n-1] \; .$$
\end{proposition}

   \noindent {\bf Preuve :} En effet, nous avons le triangle dans la cat\'egorie d\'eriv\'ee :
   $$( {\cal D}_X [s_1,\ldots, s_p] f_1^{s_1}\ldots  f_p^{s_p} )^{\ast} \rightarrow ( {\cal D}_X [s_1,\ldots, s_p] f_1^{s_1+1}\ldots  f_p^{s_p+1} )^{\ast}  
   \rightarrow  ( \frac{ {\cal D}_X [s_1,\ldots, s_p] f_1^{s_1}\ldots  f_p^{s_p}}{ {\cal D}_X [s_1,\ldots, s_p] f_1^{s_1+1}\ldots  f_p^{s_p+1}})^{\ast}[1] \; .$$
   Il r\'esulte de la proposition  pr\'ec\'edente l'isomorphisme :   
   $$ ( \frac{ {\cal D}_X [s_1,\ldots, s_p] f_1^{s_1}\ldots  f_p^{s_p}}{ {\cal D}_X [s_1,\ldots, s_p] f_1^{s_1+1}\ldots  f_p^{s_p+1}})^{\ast} \simeq
   \left(  {\cal D}_X [s_1,\ldots, s_p] f_1^{-s_1-1}\ldots  f_p^{-s_p-1} \rightarrow  {\cal D}_X [s_1,\ldots, s_p] f_1^{-s_1-2}\ldots  f_p^{-s_p-2} \right)[-n-1] \; .$$
   La proposition s'en d\'eduit.
   
          \begin{remarque} Pla\c{c}ons nous sous les hypoth\`eses $F$  localement quasi-homog\`ene 
          et   $ {\rm Der}_X (   H) $  localement libre. Le ${\cal O}_X [s_1,\ldots, s_p]$-Module  localement libre $ {\cal O}_X [s_1,\ldots, s_p] F \oplus \tilde{{\rm Der}}_X (f_1, \ldots,f_p) $
           est stable par crochet. Nous pouvons lui associer un  complexe de  Spencer que nous noterons 
           ${\rm Sp}^.( {\cal O}_X [s_1,\ldots, s_p] F \oplus \tilde{{\rm Der}}_X (f_1, \ldots,f_p)  ) $. Suivant le corollaire \ref{crr}, ce complexee   est une r\'esolution libre
           de ${\cal D}_X [s_1,\ldots, s_p]$-Module \`a gauche  
           $$   \frac{ {\cal D}_X [s_1,\ldots, s_p] f_1^{s_1}\ldots  f_p^{s_p}}{ {\cal D}_X [s_1,\ldots, s_p] f_1^{s_1+1}\ldots  f_p^{s_p+1}} \; .$$
           Son dual s'identifie naturellement au  complexe  de Spencer corespondant ${\rm Sp}^.( {\cal O}_X [s_1,\ldots, s_p] F \oplus \tilde{\tilde{{\rm Der}}}_X (f_1, \ldots,f_p)  ) $.

          \end{remarque}
          
          Soit $x_0 \in X$. Nous notons ${\cal B}(x_0, f_1, \ldots ,f_p)$ l'id\'eal de   ${\bf C}[s_1, \ldots ,s_p]$   des polyn\^omes   $b$ 
  v\'erifiant  au voisinage de $x_0$:
$$ (\ast) \;\;   b (s_1, \ldots ,s_p) \,  f_1^{s_1} \ldots  f_p^{s_p} \in {\cal D}_X[s_1, \ldots ,s_p] \,   f_1^{s_1+1} \ldots  f_p^{s_p+1} \; .$$
Nous l'appelons   id\'eal de  de Bernstein de $(f_1, \ldots ,f_p) $ au voisinage de $x_0$.\\

\begin{proposition} \label{pqhlp} Sous les hypoth\`eses $F$  localement quasi-homog\`ene  et   $ {\rm Der}_X (   H) $  localement libre .
Le  ${\cal D}_X [s_1,\ldots, s_p]$-Module \` a gauche  :
$$   \frac{ {\cal D}_X [s_1,\ldots, s_p] f_1^{s_1}\ldots  f_p^{s_p}}{ {\cal D}_X [s_1,\ldots, s_p] f_1^{s_1+1}\ldots  f_p^{s_p+1}}$$ 
est pur de nomre garde $n+1$ (tous ses sous-Modules 
ont pour nombre grade $n+1$). Pour tout $x_0$, la racine de  id\'eal de  de Bernstein de $(f_1, \ldots ,f_p) $ au voisinage de $x_0$ est principal. De plus :
$$ b(s_1, \ldots ,s_p) \in {\cal B}(x_0, f_1, \ldots ,f_p) \Longleftrightarrow b(-s_1-2, \ldots ,-s_p-2) \in {\cal B}(x_0, f_1, \ldots ,f_p)\; .$$

\end{proposition}

  \noindent {\bf Preuve :} Pour la puret\'e, suivant la remarque 5 de \cite{M}, il suffit de remarquer sous nos hypoth\`eses  :
  $$  Ext^i_{{\cal D}_X  [s_1,\ldots, s_p]}( Ext^i_{{\cal D}_X [s_1,\ldots, s_p]}
  ( \frac{ {\cal D}_X [s_1,\ldots, s_p] f_1^{s_1}\ldots  f_p^{s_p}}{ {\cal D}_X [s_1,\ldots, s_p] f_1^{s_1+1}\ldots  f_p^{s_p+1}} , {\cal D}_X [s_1,\ldots, s_p])  , {\cal D}_X [s_1,\ldots, s_p]) \neq 0 
 $$
  implique $i=n+1$. Cela  r\'esulte directement de la proposition \ref{pdq}. La proposition 20 de   \cite{M} assure que la racine de l'id\'eal de Bernstein 
  ${\cal B}(x_0, f_1, \ldots ,f_p)$ est principal. Le  r\'esultat de sym\'etrie  r\`esulte de la proposition \ref{pdq}.

\section{L'id\'eal de Bernstein d'un arrangement libre d'hyperplans lin\'eaires}

D\'esormais $V = {\bf C}^n$ ou plus g\'en\'eralement un espace vectoriel de dimension $n$. Nous consid\'erons  $H_1, \ldots ,H_p$ des hyperplans lin\'eaires
de $V$ deux \`a deux disctintcs. Ils d\'efinissent un  arrangement d'hyperplans ${\cal A} = {\cal A} ( H_1, \ldots ,H_p  ) $ de $V$.
Notons toujours $H$ la r\'eunions de ces hyperplans. Tous supposerons l'arrangement d\'efini par $H$ libre.
Pour $1\leq i\leq n$, nous consid\'erons $l_i$ une forme lin\'eaire homog\`ene sur $V$ de noyaux $H_i$, posons $L = l_1 \ldots l_p$. Les formes lin\'maires $l_i$ sont d\'efinies \`a un coefficient de proportionalit\'e
 pr\`es. Le ${\cal O}_X $-Module,   $ {\rm D} (  {\cal A} ) = {\rm Der}_{V}  (    L)$ est donc un ${\cal O}_V$-Module libre.\\
 
 Suivons les notations du livre de P. Orlik et H. Terao \cite{O-T}. Notons $L( {\cal A} )$  l'ensemble des intersections non vide  d'\'el\'ements de  $ \{ H_1, \ldots ,H_p \} $.
  Si $X \in L( {\cal A} )$, nous notons   $ J(X) =\{ i \in \{1, \ldots ,p\} \; ; \; X \subset H_i\} $,    $r(x)$ le codimension de $X$  
  et  ${\cal A}_X $ l'arrangement d'hyperplans d\'efini par $ \{ H_i \}_{ i \in J(X) } $. Suivant \cite {O-T}, si ${\cal A}$  est libre, l'arrangement
 ${\cal A}_X $ reste libre. \\

Notons $E$ le champ d'Euler. Dans un syst\`eme de coordonn\'ees $(x_1,\ldots , x_n)$ de $X$ :
  $ E = \sum_{j=1}^n x_j \frac{\partial }{\partial x_j}$  Le champ d'Euler est un champ de vecteur logarithmique  : $E(L)=pL$ et $E(l_i)= l_i$. Suivant les notations du paragraphe pr\'ec\'edent,
  nous avons $\tilde{E} = E - \sum_{j=1}^p s_j$.\\
  
  Une d\'erivation $\chi$
   est homog\`ene de degr\'e $r$, si elle  s'\'ecrit :
   $$ \chi= \sum_{i=1}^n a_i(x) \frac{\partial }{\partial x_i}\;  ,$$
   o\`u les $a_i(x) $ sont des polyn\^omes homog\`enes de degr\'e $r$.\\
   
   Suivant \cite{O-T},  ${\rm Der}_{{\cal O}_V} ( L)$ admet une base de champs de vecteurs $( \delta _1, \ldots , \delta _n )$  homog\`enes. De plus, la suite
   $({\rm deg} (\delta _1), \ldots , {\rm deg} (\delta _n))$ des degr\'es d'homog\'en\'eit\'e des $\delta _i$  est ind\'ependant de la base
   et est appel\'ee la suite des exposants de l'arrangement   $ {\cal A}$. Nous   notons ${\rm exp} ( {\cal A})  = ( 0^{e_0}, 1^{e_1}, 2^{e_2},  \ldots ) $ pour signifier que dans une
   base, il y a $e_i$  d\'erivations  homog\`enes de degr\'e $i$. De plus, nous avons $e_0=0$ si et seulement si le rang de la famille de formes lin\'eaires
   $(l_1,  \ldots ,l_p)$ est \'egal \`a $n$.  Si $e_0 >0$, l'arrangement  $  {\cal A}  $ est produit de l'arrangement $({\bf C}^{e_0}, \emptyset) $ et d'un arrangement libre
   $({\bf C}^{n-e_0},  {\cal A}') $ o\`u  ${\rm exp} ( {\cal A}')  = ( 0^{0}, 1^{e_1}, 2^{e_2},  \ldots ) $.  Si $e_0$ est nul, l'arrangement  $ {\cal A} $ est isomorphe au produit  de $e_1$
   arrangements libres irr\'eductibles de rang maximal. L'arrangement $ {\cal A}$ est dit irr\'eductible s'il n'est pas isomorphe \` a un produit d'arrangements   non vides. 
   Une caract\'erisation d'un arrangement irr\'eductible est donc que son   exposant $e_1$ soit \'egal \`a $1$. \\

Soit $A_n({\bf C})$, l'alg\`ebre de Weyl des op\'erateurs diff\'erentiels \`a coefficients polynomiaux sur  $V$. Suivant la d\'emonstration de Bernstein \cite{Be}, l'id\'eal des polyn\^omes
 $ b\in {\bf C}[s_1, \ldots ,s_p]$ v\'erifiant :
 $$   \;\;   b (s_1, \ldots ,s_p) \,   l_1^{s_1} \ldots  l_p^{s_p} \in A_n({\bf C})[s_1, \ldots ,s_p] \,   l_1^{s_1+1} \ldots  l_p^{s_p+1} \; ,$$
 n'est pas  r\'eduit \`a z\'ero. Il ne d\'epend pas du choix des formes lin\'eaires  $l_i$ qui d\'efinissent les hypersurfaces $H_i$. Nous notons cet id\'eal 
ou  ${\cal B} ({\cal A} ( H_1, \ldots ,H_p  ))$ et l'appelons l'id\'eal de Bernstein
 de ${\cal A} $.  \\
 
 Pour des raisons d'homog\'en\'eit\'e : ${\cal B} ({\cal A} ( H_1, \ldots ,H_p  )) = {\cal B} ( 0, l_1, \ldots ,l_p  ) $.  De plus, si $( \delta _1, \ldots , \delta _n )$  
 est une base  de champs de vecteurs    homog\`enes de  ${\rm Der}_{{\cal O}_V} ( L)$, nous avons montre \`a la propsoition  \ref{pqhlp} :
$$ {\rm Ann}_{ A_n({\bf C}) [s_1,\ldots, s_p]}  l_1^{s_1}\ldots  l_p^{s_p}
    = {\cal D}_V [s_1,\ldots, s_p] (   \tilde{\delta}_1 ,\ldots ,     \tilde{\delta}_l )  \; . $$ 
    De plus, l'arrangement ${\cal A}$ est localement quasi-homog\'en\-e. Nous avons donc suivant le paragraphe pr\'ec\'edent : 
 $$   \frac{ {\cal D}_V [s_1,\ldots, s_p] l_1^{s_1}\ldots  l_p^{s_p}}{ {\cal D}_V [s_1,\ldots, s_p] l_1^{s_1+1}\ldots  l_p^{s_p+1}}$$ 
est pur de nombre grade $n+1$ (tous ses sous-Modules 
ont pour nombre grade $n+1$) et la   racine de  id\'eal    ${\cal B} ({\cal A} ( H_1, \ldots ,H_p  )) $ est principal. De plus, nous avons la propri\'et\'e de 
sym\'etrie :
$$ b(s_1, \ldots ,s_p) \in {\cal B} ({\cal A} ( H_1, \ldots ,H_p  ))   \Longleftrightarrow b(-s_1-2, \ldots ,-s_p-2) \in{\cal B} ({\cal A} ( H_1, \ldots ,H_p  ))  \; .$$

\begin{lemme} \label{lde} Supposons que ${\cal A} ( H_1, \ldots ,H_p  ) $ soit un arrangement libre  et  irr\'eductible. Soit $ r  $ la codimension  de $H_1 \cap \ldots \cap H_p$, 
alors tout $b \in  {\cal B} ({\cal A} ( H_1, \ldots ,H_p  ))$ est multiple de
$$  \prod_{j=0}^{2(p - r)}(s_1 + \cdots +s_p + n +j )\;  .$$
\end{lemme}

  \noindent {\bf Preuve :}  Comme l'arrangement est suppos\'e irr\'educible $e_1=1$. Commencons par supposer $r=n$ ou  encore  $e_0 = 0$. \\
  
  Soit $r = {\rm inf} \, \{ i  \in  {\rm exp} ( {\cal A}) \; ; \; i \geq 2\}$. Montrons  tout d'abord  que 
$\prod_{j=0}^{r+p-3}(s_1 + \cdots +s_p + n +j )$ divise tout $b$ de l'id\'eal de Bernstein.  Soit $v$ un polyn\^ome homog\`ene de degr\'e $0\leq d \leq r-2$. Si $b \in {\cal B} ({\cal A} ( H_1, \ldots ,H_p  )) $, en multipliant l'identit\'e fonctionnelle
associ\'ee \` a $b$ par $v l_1 \ldots l_k$ o\`u $k\leq p-1$, nous obtenons en particulier  :
$$  b(s_1, \ldots ,s_p) v l_1^{s_1+1} \ldots  l_k^{s_k+1}l_{k+1}^{s_{k+1}} \ldots  l_p^{s_p} \in A_n({\bf C})[s_1, \ldots ,s_p] \, l_1^{s_1+1} \ldots  l_k^{s_k+1}l_{k+1}^{s_{k+1}+1 } l_{k+2}^{s_{k+2} } \ldots  l_p^{s_p}   \; ,$$ 
Donc, il existe un op\'orateur diff\'erentiel $ P \in A_n({\bf C}) [s_1,\ldots, s_p]$ tel que
$$  b(s_1, \ldots ,s_p) v - Pl_{k+1} \in  {\rm Ann}_{ A_n({\bf C}) [s_1,\ldots, s_p]}  l_1^{s_1+1} \ldots  l_k^{s_k+1}l_{k+1}^{s_{k+1}} \ldots  l_p^{s_p} \; .$$
Ainsi, il existe des $ A_i \in A_n({\bf C}) [s_1,\ldots, s_p]$ tel que :
$$(\ast)\;   b(s_1, \ldots ,s_p) v  =   Pl_{k+1} + A_1 \tilde{\delta  _1}(s_1 +1, \ldots, s_k +1 , s_{k+1},\ldots s_p) + \cdots + A_n \tilde{\delta  _n }(s_1 +1, \ldots, s_k +1 , s_{k+1},\ldots s_p) \; .$$

Un op\'erateur de ${\cal D}_{ {\bf C}^n ,0 }  [s_1, \ldots ,s_p]$ (resp  $ A_n({\bf C}) [s_1,\ldots, s_p]$) s'\'ecrit de fa\c{c}on unique :
 $$P = \sum_{\alpha \in A } \frac{\partial}{\partial x^{\alpha} }c_{\alpha }(s_1, \ldots ,s_p) \; , $$
 o\`u $c_{\alpha  }(s_1, \ldots ,s_p)  \in {\cal O}_{ {\bf C}^n  ,0 }  [s_1, \ldots ,s_p] $ (resp.   ${\bf C}[x_1, \ldots, x_n,s_1, \ldots ,s_p]$)  et $A$ est un ensemble fini de ${\bf N}^n$.
 Cette \'ecriture sera dite l'\'ecriture \`a droite de $P$.
 Nous appelons $c_{0}(s_1, \ldots ,s_p) $ le terme constant de l'\'ecriture \`a droite de $P$.\\

Prenons la partie homog\`ene de degr\'e $d$  des deux  termes constants de l'\'ecriture \`a droite  dans l'\'egalit\'e $(\ast)$.
Sachant que l'on peut supposer  que $\delta  _1$ est le champ d'Euler $E$, nous obtenons  :
$$   b(s_1, \ldots ,s_p) v  = p  l_{k+1} + (\sum_{j=1}^p s_j +k +d +n)a_1 \;l ,$$
o\`u $p$ et $a_1$ sont des polyn\^omes homog\`enes de degr\'e $d-1$ (resp. $d$) et plus pr\'ecisement les parties homog\`unes de degr\'e $d-1$ (resp. $d$) des termes constants \`a droite   de $P$ et $A_1$.
En choisissant $v$ non multiple
 de $l_{k+1}$, nous obtenons : $\sum_{j=1}^p s_j +k +d +n$ divise $b$.   Ainsi, $\prod_{j=0}^{r+p-3}(s_1 + \cdots +s_p + n +j )$ divise 
 $b$. Or $ n+p + r -3 \geq p-n$, il reste \`a utiliser  la propri\'et\'e de sym\'etrie de  ${\cal A} ( H_1, \ldots ,H_p  ) $.\\
 
 Si $e_0> 0$, $r = n -e_0$ : L'arrangement est le produit de l'arrangement vide sur $ {\bf C}^{e_0}$  et d'un arrangement $(H'_1, \ldots , H'_p)$ de $ {\bf C}^{n- e_0}$   d'exposant $(0,1,e_2, \ldots, )$. Il est facile de montrer que
 l'id\'eal de Bernstein de l'arrangement ${\cal A} ( H_1, \ldots ,H_p  ) $ co\"{\i}ncide avec celui de l'arrangement  ${\cal A} (H'_1, \ldots , H'_p) $. 
 Il reste \`a utiliser  le premier cas $r=n$ que nous venons de traiter.

 \begin{proposition} (diviseur \'evident) Supposons que ${\cal A} ( H_1, \ldots ,H_p  ) $ soit un arrangement libre, tout  polyn\^ome de Bernstein  est multiple de  :
 $$b_{de}(s_1, \ldots ,s_p) = \prod_{X \in   L( {\cal A} ) \; ; \; {\cal A}_X \; irr\acute{e}ductible}   \prod_{ j=0}^{ 2 ({\rm card} J(X) - r(X) )}  ( \sum_{j \in J(X) }  s_j  + r(X) + j) \; .$$ 
 \end{proposition}

    \noindent {\bf Preuve :} Soit $X \in   L( {\cal A} )$ tel que ${\cal A}_X$ soit irr\'eductible. Nous savons voir \cite{O-T} que ${\cal A}_X$ reste un arrangement libre.
 Pour $j \notin J(X) $, $l_j$ est non identiquement nul sur $X$. Nous pouvons trouver $x_0 \in X$
    tel que pour  $j \notin J(X) $, $f_j (x_0) \neq 0$ et donc tel que  le germe $f_j$ en $x_0$ soit inversible. Il alors facile de montrer que l'id\'eal de Bernstein ${\cal B}(x_0, f_1, \ldots ,f_p)$
    est engendr\'e par les  polyn\^omes de  l'id\'eal de Bernstein ${\cal B}(x_0, (f_j)_{j \in  J(X)})$. Tous ces polynomes sont d'apr\`es le lemme \ref{lde} multiples de : 
   $$ \prod_{ j=0}^{ 2 ({\rm card} J(X) - r(X) )}  ( \sum_{j \in J(X) }  s_j  + r(X) + j) \; ,$$ 
   ce qui d\'emontre la proposition.

    \begin{proposition} (multiple \'evident) Soit  ${\cal A} ( H_1, \ldots ,H_p  ) $   un arrangement d'hyperplan lin\'eaire {\bf non n\'ecessairement libre}, il existe un entier $N$ assez grand tel que :
 $$b_{me}(s_1, \ldots ,s_p) = \prod_{X \in   L( {\cal A} ) \; ; \; {\cal A}_X \; irr\acute{e}ductible}   \prod_{j=0}^{ N} ( \sum_{j \in J(X) }  s_j  + r(X) + j) \; .$$ 
 soit un polyn\^ome de Bernstein de  ${\cal A} ( H_1, \ldots ,H_p  ) $.
 \end{proposition}

  \noindent {\bf Preuve :} Nous pouvons supposer l'arrangement $ {\cal A} (H_1, \ldots ,H_p ) $ irr\'eductible et  que l'intersection des hypersurfaces de  l'arrangement soit r\'eduit au vecteur nul.
 Ainsi,   $r(H_1 \cap \ldots \cap H_p ) =n$.\\
 
      Soit $i_1,\ldots ,i_n \in {\bf N}$ tels que $i_1+ \cdots + i_n = k$. Notons $C_k^{i_1,\ldots ,i_n}$ la suite d'entiers 
      d\'efinie par 
      r\'ecurrence par :
      $$ C_k^{i_1,\ldots ,i_n}= C_{k-1}^{i_1-1,\ldots ,i_n}+\cdots + C_{k-1}^{i_1 ,\ldots ,i_n-1} \; , \;  C_1^{0,\ldots, 01,0,\ldots, 0}  =1    \; .$$
      
      D\'esignons par $(x_1, \ldots , x_n) $ le syst\`eme de coordon\'ees canoniques de $V$. Nous v\'erifions par r\'ecurrence  :
      $$ \prod_{j=0}^{k-1} ( \sum_{i=1}^n x_i \frac{\partial}{\partial x_i}   - j) = \sum_{i_1+ \cdots + i_n = k}  C_k^{i_1,\ldots ,i_n} x_1^{i_1} \cdots  x_n^{i_n} 
  \frac{\partial ^{i_1}}{\partial x_1^{i_1}}   
     \cdots  \frac{\partial ^{i_n}}{\partial x_n^{i_n}} \; .$$
      D'o\`u par transposition :
     $$ \prod_{j=0}^{k-1} ( \sum_{i=1}^n x_i \frac{\partial}{\partial x_i}   + j + n) 
      = \sum_{i_1+ \cdots + i_n = k}  C_k^{i_1,\ldots ,i_n}   
  \frac{\partial ^{i_1}}{\partial x_1^{i_1}}   
     \cdots  \frac{\partial ^{i_n}}{\partial x_n^{i_n}}  x_1^{i_1} \cdots  x_n^{i_n} \; .$$
     
      Il en r\'esulte :
\begin{lemme}\label{lmp} Si la famille de formes lin\'eaires $(l_1, \ldots ,l_p)$ est de rang $n$, 
nous avons :  :
 $$\prod_{j=0}^{k-1} ( s_1+\cdots +s_p +n+j) \; l_1^{s_1} \ldots  l_p^{s_p}     \in  
 ( \sum_{i_1+ \cdots + i_n = k}  C_k^{i_1,\ldots ,i_n}   
  \frac{\partial ^{i_1}}{\partial x_1^{i_1}}   
     \cdots  \frac{\partial ^{i_n}}{\partial x_n^{i_n}} ) ( l_1, \ldots ,l_p)^k   \; l_1^{s_1} \ldots  l_p^{s_p} \; .$$
      \end{lemme}
  Notons  $ L''( {\cal A} ( H_1, \ldots ,H_p    )$ l'ensemble des  $X \in L( {\cal A} )$ tels que  $r(X)  = n-1$. Nous avons par r\'ecurrence :
  $$ ( l_1, \ldots ,l_p)^{p-n+1} \in  (\prod_{j \notin J(X)}l_j)_{X \in L''( {\cal A} ( H_1, \ldots ,H_p    )  }  \; .$$
 Il en r\'esulte que pour tout $k$ il existe un entier $K$ assez grand tel que :
 $$\prod_{j=0}^{K-1} ( s_1+\cdots +s_p +n+j) \; l_1^{s_1} \ldots  l_p^{s_p}     \in  
 ( \sum_{i_1+ \cdots + i_n = K}  C_k^{i_1,\ldots ,i_n}   
  \frac{\partial ^{i_1}}{\partial x_1^{i_1}}   
     \cdots  \frac{\partial ^{i_n}}{\partial x_n^{i_n}} )(\prod_{j \notin J(X)}l_j)^k_{X \in L''( {\cal A} ( H_1, \ldots ,H_p    )}    \; l_1^{s_1} \ldots  l_p^{s_p} \; .$$
 
 Par r\'ecurrence sur  la dimension, pour tout $X \in L''( {\cal A})$, il existe $N_X$ entier assez grand tel  que :
 $$\prod_{Y \in   L( {\cal A}_X ) \; ; \; {\cal A}_Y \; irr\acute{e}ductible}   \prod_{j=0}^{ N_X} ( \sum_{j \in J(Y) }  s_j  + r(Y) + j)
 \prod_{j \in J(X)}l_j^{s_j}
      \in  {\cal D}_X [s_1,\ldots, s_p] \prod_{j \in J(X)}l_j^{s_j+1}\; .$$
      Donc, si  $k$ est un entier assez grand :
      $$    \prod_{Y \in   L( {\cal A}_X ) \; ; \; {\cal A}_Y \; irr\acute{e}ductible}   \prod_{j=0}^{ N_X} ( \sum_{j \in J(Y) }  s_j  + r(Y) + j)
  (\prod_{j \notin J(X)}l_j)^k    l_1^{s_1} \ldots  l_p^{s_p} \in A_n({\bf C})[s_1, \ldots ,s_p] \,   l_1^{s_1+1} \ldots  l_p^{s_p+1} \; ,$$   
  La proposition s'en d\'eduit.

\begin{theoreme}  Supposons que ${\cal A} ( H_1, \ldots ,H_p  ) $ soit un arrangement libre  d'hyperplans lin\'eaires. Alors,  l'id\'eal de Bernstein
 ${\cal B} ({\cal A} ( H_1, \ldots ,H_p  ))$ est principal de g\'en\'erateur :
 $$  \prod_{X \in   ) \; ; \; {\cal A}_X \; irr\acute{e}ductible}   \prod_{ j=0}^{ 2 ({\rm card} J(X) - r(X) )}  ( \sum_{j \in J(X) }  s_j  + r(X) + j) \; .$$
 \end{theoreme}

     \noindent {\bf Preuve :} Suivant la proposition    \ref{pqhlp}, la racine de ${\cal B} ({\cal A} ( H_1, \ldots ,H_p  ))$ est principal. 
     Consid\'erons un g\'en\'erateur $c$ de cet id\'eal.
     l'id\'eal  ${\cal B} ({\cal A} ( H_1, \ldots ,H_p  ))$
     v\'erifie la propri\'et\'e de sym\'etrie :
      $$ b(s_1, \ldots ,s_p) \in {\cal B} ({\cal A} ( H_1, \ldots ,H_p  )) \Longleftrightarrow b(-s_1-2, \ldots ,-s_p-2) \in {\cal B} ({\cal A} ( H_1, \ldots ,H_p  ))\; .$$
      Sa racine v\'erifie donc cette propri\'et\'e. Il en r\'esulte que 
      $  c(-s_1-2, \ldots ,-s_p-2)$ est \'egal \`a plus ou moins    $c(s_1, \ldots ,s_p)$.\\
      
       Puisque $b_{de}$ est un polyn\^ome r\'eduit, le polynome $c$ est clairement  multiple du  polyn\^ome   $b_{de}$. 
     Il divise d'autre part $b_{me}$.  Supposons alors que $\sum_{j \in J(X) }  s_j  + r(X) + j$ avec 
     $j>  2 ({\rm card} J(X) - r(X) )$ soit un facteur de $c$.  Le polyn\^ome $c$ aurait donc comme facteur :\\
     $$- \sum_{j \in J(X) }  s_j - 2 {\rm card} J(X)  + r(X) + j= - ( \sum_{j \in J(X) }  s_j  + r(X)  +  2 ({\rm card} J(X) - r(X) ) - j\; .$$
     C'est impossible puisque comme $r(X)  +  2 ({\rm card} J(X) - r(X) ) - j<0$ et que  ce facteur n'apparait pas dans  $b_{me}$.  Il en r\'esulte que la  racine de ${\cal B} ({\cal A} ( H_1, \ldots ,H_p  ))$ est principal
     de g\'en\'erateur $b_{de}$.\\
     
     Ecrivons  $b_{me}(s_1, \ldots ,s_p) = b_{de}(s_1, \ldots ,s_p) a(s_1, \ldots ,s_p)$.  Regardons le $ {\cal D}_X [s_1,\ldots, s_p] $ sous-Module :

  $$ L' = b_{de}(s_1, \ldots ,s_p) \frac{ {\cal D}_X [s_1,\ldots, s_p] f_1^{s_1}\ldots  f_p^{s_p}}{ {\cal D}_X [s_1,\ldots, s_p] f_1^{s_1+1}\ldots  f_p^{s_p+1}}  \subset 
L=  \frac{ {\cal D}_X [s_1,\ldots, s_p] f_1^{s_1}\ldots  f_p^{s_p}}{ {\cal D}_X [s_1,\ldots, s_p] f_1^{s_1+1}\ldots  f_p^{s_p+1}} /; .$$
Ce ${\cal D}_X [s_1,\ldots, s_p] $-Module est support\'e par les z\'eros du  polyn\^ome $ a(s_1, \ldots ,s_p)$. Nous venons de voir que  $L$ est support\'e par les z\'eros du  polyn\^ome   $ b_{de}(s_1, \ldots ,s_p$.
Il en r\'esulte que la dimension de la vari\'et\'e caract\'eristique de $L'$ est strictment inf\'erieur \`a celle de $L$. Cela contredit la puret\'e de $L$
\'etablie \` a la proposition  \ref{pqhlp}.\\

Terminons par un r\'esultat  sur les composantes irr\'eductibles de $W^{\sharp}_{l_1 , \ldots l_p } \cap H$. Soit $X\in  L( {\cal A})$. Notons pour $V= {\bf C}^n$ :
$$\Omega _{ T^{\ast}_X {\bf C}^n , (l_j)_{ j \notin J(X)} } = \{ ( x, \xi +  \sum_{j \notin J(X) }     s_j \frac{dl_j}{l_j} , (s_j)_{j \notin J(X) }  ) \; ;\; (x,\xi) \in   T^{\ast}_X {\bf C}^n \; {\rm et} \; 
 \prod_{j \notin J(X)} l_j(x) \neq 0 \} \; .$$
 Notons  $W^{\sharp}_{X, (f_j)_{ j \notin J(X)} }$ son adh\'erence. C'est un ensemble irr\'eductible de dimension de $ T^{\ast}  {\bf C}^n \times {\bf C}^{p- {\rm card} \, J(X)}$ de dimension
 $n+p - {\rm card} J(X)$ qui s'identifie \`a $W^{\sharp}_{{f_1}_{\mid X} , \ldots , {f_p}_{\mid X} }\times {\bf C}^{n -  r(X) }$.\\

Sans aucune hypoth\'ese  suivant \cite{B-M-M1}, les composantes irr\'eductibles de $W^{\sharp}_{f_1 , \ldots f_p } \cap H$
  sont de dimension $n+p-1$ et se projette apr la projection $(x, \xi,s) \mapsto s$   sur   des hyperplans lin\'eaires appel\'es    pentes de  $ f_1 , \ldots f_p$.\\

\begin{proposition} Supposons que ${\cal A} ( H_1, \ldots ,H_p  ) $ soit un arrangement libre  d'hyperplans lin\'eaires.
Les composantes irr\'eductibles de  $W^{\sharp}_{l_1 , \ldots l_p } \cap H$ sont les :
$$   W^{\sharp}_{X, (l_j)_{ j \notin J(X)} } \times \{ (s_j)_{ j \in J(X) }   \in {\bf C}^{{\rm card} J(X)} \; ; \; \sum_{j \in J(X)} s_j = 0  \}$$
o\`u $X$ parcourt les \'el\'ements de $L( {\cal A}) $ tels que ${\cal A}_X$ soit irr\'eductible.
 \end{proposition}
 
      \noindent {\bf Preuve :}  Soit $X_0 = H_1 \cap \ldots \cap H_p$.  Soit $(e_0,e_1, ...)$ les exposants de l'arrangement ${\cal A}$.   Si $e_1 >1$, 
    aucune composante irr\'eductible de   $W^{\sharp}_{l_1 , \ldots , l_p } \cap H$  n'est contenue dans $X_0$. Si $e_1=1$ :   
    $$W^{\sharp}_{l_1 , \ldots ,l_p } \cap X_0 =  T^{\ast}_{X_0} {\bf C}^n \times  \{(s_;, \ldots ,s_p) \in  {\bf C}^p \; ; \; \sum_{j=1}^p s_j=0 \} \; .$$
   Pour d\'emontrer cela \`a un facteur trivial pr\`es, nous pouvons supposer $e_0=0$. Le r\'esultat provient alors du fait que $W^{\sharp}_{l_1 , \ldots ,l_p }$ est d\'efini par l'id\'eal
  $  ( \sigma^{\sharp} 
     ( \tilde{\delta}_1 ),\ldots ,   \sigma^{\sharp}   ( \tilde{\delta}_n ) )$ o\`u   $  (  {\delta}_1 ,\ldots ,    {\delta}_n)$ est une base de 
      ${\rm Der}_{ {\cal O}_{ {\bf C}^n } }( L)$.
  et que   le champ d'Euler peut  \^etre choisi dans cette famille de g\'en\'erateurs. Dans ce cas,  $W^{\sharp}_{l_1 , \ldots , l_p } \cap H$
   est de dimension $n+p-1$. C'est donc une composante irr\'eductible de  $W^{\sharp}_{l_1 , \ldots l_p } \cap H$.\\
   
 Soit $X \in L( {\cal A}) $. Cherchons les composantes irr\'eductibles de  $W^{\sharp}_{l_1 , \ldots l_p } \cap H$ contenues dans $X$ et non contenues
 dans $X'\subset X$ et $X' \in L( {\cal A}) $. Un telle composante est contenue dans l'adh\'erence  de :
 $$ \{ (x, \xi + \sum_{j \notin J(X)} s_j \frac{dl_j}{l_j}, s_1, \ldots ,s_p)  \; ; \; (x,\xi, (s_l)_{j \in J(X)} ) \in   W^{\sharp}_{(l_j)_{j \in J(X)  } } \cap X
 \} \; . $$
 D'apr\`es le cas pr\'ec\'edent,  si  ${\cal A}_X$ n'est pas  irr\'eductible  cette adh\'erence n'a pas la bonne dimension  . Et si   ${\cal A}_X$  est    irr\'eductible, cette adh\'erence s'identifie \`a :
 $$   W^{\sharp}_{X, (l_j)_{ j \notin J(X)} } \times \{ (s_j)_{ j \in J(X) }   \in {\bf C}^{{\rm card} J(X)} \; ; \; \sum_{j \in J(X)} s_j = 0  \}$$
 qui est de dimension $n+p-1$. C'est donc une composante irr\'eductible de $W^{\sharp}_{l_1 , \ldots l_p } \cap H$.
 
\begin{remarque}
 Nous retrouvons que si ${\cal A} ( H_1, \ldots ,H_p  ) $ est un arrangement libre  d'hyperplans lin\'eaires, les pentes de  $ l_1 , \ldots l_p$ sont les hyperplans de 
 ${\bf C}^p$ d'\'equations $ \sum_{j \in J(X)} s_j = 0  $ o\`u $X$ parcourt les \'el\'ements de $L( {\cal A}) $ tels que ${\cal A}_X$ soit irr\'eductible.
\end{remarque}

 \end{document}